\documentclass{amsart}
\usepackage{amssymb}
\newtheorem{thm}{Theorem}
\newtheorem{lem}{Lemma}

\begin{document}

\bibliographystyle{plain}

\title[Romi F. Shamoyan and Milo\v s Arsenovi\'c]{Some remarks on extremal problems in weighted Bergman spaces of analytic functions}

\author[]{Romi F. Shamoyan}
\author[]{Milo\v s Arsenovi\' c$\dagger$}

\address{Bryansk University, Bryansk Russia}
\email{\rm rshamoyan@yahoo.com}

\address{Faculty of mathematics, University of Belgrade, Studentski Trg 16, 11000 Belgrade, Serbia}
\email{\rm arsenovic@matf.bg.ac.rs}

\thanks{$\dagger$ Supported by Ministry of Science, Serbia, project M144010}

\date{}

\maketitle

\begin{abstract}

We prove some sharp extremal distance results for functions in weighted Bergman spaces on the upper halfplane.
We also prove new analogous results in the context of bounded strictly pseudoconvex domains with smooth boundary.

\end{abstract}

\footnotetext[1]{Mathematics Subject Classification 2010 Primary 30D45.  Key words
and Phrases: Bergman spaces, pseudoconvex domains, extremal problems.}

\section{Introduction}

If $Y$ is a normed space and $X \subset Y$, then we set ${\rm dist}_Y (f, X) = \inf_{g \in X} \| f - g \|_Y$. If the
space $Y$ is clear from the context, we write simply ${\rm dist}(f, X)$. In the problems we are going to consider,
$X$ itself is going to be a (quasi)-Banach space.

We denote by $H(\Omega)$ the space of all holomorphic functions on an open set $\Omega \subset \mathbb C^n$.
In this paper we consider distance problems in weighted Bergman spaces over the upper half-plane
$\mathbb C_+ = \{x + iy : y > 0 \}$ and over a bounded strictly pseudoconvex domain $\Omega \subset \mathbb C^n$
with smooth boundary.

The weighted Bergman space $A^p_\alpha (\mathbb C_+)$ consists of all functions $f\in H(\mathbb C_+)$ such that
$$\| f \|_{p,\alpha} = \left(\int_0^\infty \int_{-\infty}^\infty |f(x+iy)|^p y^\alpha dx dy \right)^{\frac{1}{p}}
< \infty,$$
where $\alpha > -1$ and $0 < p < \infty$ (see \cite{MDAD} and \cite{DK}). The above spaces are Banach spaces for $p \geq 1$ and complete metric spaces for $0<p<1$. It is natural to consider the space $A^\infty_\nu = A^\infty_\nu(\mathbb C_+)$ of all holomorphic functions $F$ on $\mathbb C_+$ satisfying
$$\| F \|_{A^\infty_\nu} = \sup_{y>0} \sup_{x \in \mathbb R} |F(x+iy)|y^{\nu} < \infty,$$
where $\nu > 0$, this space is also a Banach space.

Let $D$ be a bounded strictly pseudoconvex domain in $\mathbb C^n$ with smooth boundary and let $\delta(z)$ denote the distance from $z \in D$ to the boundary of $D$ with respect to some Riemannian metric, see \cite{Be}. We define
$$L^p_s(D) = L^p(D, \delta(z)^{s-\frac{n+1}{p}}dV(z)),$$
where $sp>n$, $0<p<\infty$ and $dV(z)$ is the volume element on $D$, (see \cite{Be}). Set
$$A_s^p(D) = L^p_s(D) \cap H(D),\;\;\;\; sp>n, \;\; 0<p<\infty,$$

$$A^\infty_s(D) = \{ f \in H(D) : \sup_{z\in D} |f(z)|\delta^s (z) < \infty \},\;\; s \geq 0.$$
These spaces are Banach spaces for $1 \leq p \leq \infty$ and complete metric spaces for $0<p<1$.

In this paper we investigate the following two problems: estimate the distance, in
$A^\infty_{\frac{\nu+2}{q}}(\mathbb C_+)$ norm,
from $f \in A^\infty_{\frac{\nu+2}{q}}(\mathbb C_+)$ to $A^q_\nu(\mathbb C_+)$ and estimate the distance,
in $A^\infty_s(D)$ norm, from $f \in A^\infty_s(D)$ to $A^q_s$. In both cases we give sharp results. The next section
deals with the upper half-plane case, the last one deals with strongly pseudoconvex domains. Techniques used to obtain
our results in these two different settings are similar, and the same ideas were used to study analogous problems for analytic Besov spaces in the unit ball and polydisc (see \cite{ShM} and \cite{SM}). The literature on the extremal
problems in spaces of analytic functions is extensive, even in the case of the unit disk, a classical exposition of
these problems treated by duality methods developed by S. Havinson, W. Rogosinski and H. Shapiro can be found in \cite{Du}.

\section{Distance problems in $A^p_{\frac{\nu + 2}{p}}(\mathbb C_+)$ spaces}

The main tool in our investigation is Bergman representation formula. We first collect the results needed in the proofs
of our theorems. The following result is contained in \cite{MDAD}:
\begin{lem}\label{linre}
If $f \in A^p_\alpha (\mathbb C_+)$, $0<p<\infty$ and $\alpha > -1$ then
$$f(z) = \frac{\beta + 1}{\pi} \int_0^\infty \int_{-\infty}^\infty
\frac{f(w) (\Im w)^\beta}{(\overline w - z)^{2+\beta}} dm_2(w),$$
where $0 < p \leq 1$, $\beta \geq \frac{2+\alpha}{p} -2$ or $1 \leq p < \infty$, $\beta \geq
\frac{1+\alpha}{p}-1$.
\end{lem}
Here and in the following $m_2$ denotes two dimensional Lebesgue measure. The next two results are contained in \cite{DK}.

\begin{lem}\label{lcub}
Let $f \in H(\mathbb C_+)$. Assume $\mathbb C_+$ is covered by dyadic squares: $\mathbb C_+ = \cup_{k=1}^\infty
\Delta_k$ and let $(\Delta_k^\star)$ be the corresponding family of enlarged squares (see \cite{DK}). Then
$\Delta_k^\star$ is a finitely overlapping covering of $\mathbb C_+$ and
$$\sup_{z \in \Delta_k} |f(z)|^p (\Im z)^\alpha \leq \frac{C}{|\Delta_k^\star|} \int_{\Delta_k^\star}
|f(z)|^p (\Im w)^\alpha dm_2(w)$$
for $0<p<\infty$ and $\alpha > -1$.
\end{lem}

\begin{lem}\label{lines}
Let $\Delta_k$ and $\Delta_k^\star$ are as in the previous lemma, let $w_k$ be the center of the diadic square $\Delta_k$. Then we have:
$$(\Im w_k)^2 \asymp | \Delta_k | = m_2(\Delta_k) \asymp m_2(\Delta_k^\star),$$

$$|\overline w - z | \asymp |\overline w_k - z|, \;\; w\in \Delta_k \;\;\;, z \in \mathbb C_+,$$

$$\Im w \asymp \Im w_k, \;\;\; w \in \Delta_k$$
and the following integral estimate:
$$\int_{\mathbb C_+} \frac{(\Im z)^\alpha dm_2(z)}{|\overline w - z|^{(2+\beta)p}} \leq C
(\Im w)^{\alpha + 2 - (\beta + 2)p},\;\;\; w \in \mathbb C_+$$
valid for all $\beta$ satisfying $(\beta + 2)p - 2 > \alpha$, $\alpha > -1$.
\end{lem}

\begin{thm} {\rm (see \cite{BB} and \cite{BBG})} Let $0< p < \infty$, $\nu > 0$. Then there is a constant
$C = C(p, \nu) > 0$ such that for all $x+iy \in \mathbb C_+$ and all $F \in A^p_\nu(\mathbb C_+)$ we have
$$|F(x+iy)| \leq C y^{-\frac{\nu + 2}{p}} \| F \|_{A^p_\nu(\mathbb C_+)}.$$
Also, there is a constant $C = C(p, \nu) > 0$ such that for all $y > 0$ and all $F \in A^p_\nu(\mathbb C_+)$ we have
$$\left(\int_{-\infty}^\infty |F(x+iy)|^p dx \right)^{1/p} \leq C y^{-\frac{\nu + 1}{p}} \| F \|_{A^p_\nu(\mathbb C_+)}.$$
\end{thm}

We clearly have
$$\| F \|_{A^\infty_{\frac{\nu + 2}{p}}} \leq C \| F \|_{A^p_\nu}$$
for $0<p<\infty$, $\nu > 0$ and $F \in A^p_\nu(\mathbb C_+)$.

The next two theorems show that ${\rm dist}_{A^q_{\frac{\nu + 2}{q}}}(f, A^q_\nu)$ can be explicitly given in the case
of the upper halfplane $\mathbb C_+$. We treat separately cases $0<q\leq 1$ and $q>1$. The distance we are looking for
is described using the following sets:
$$V_{\epsilon, t}(f) = \{ z = x+iy \in \mathbb C_+ : |f(x+iy)|y^t \geq \epsilon \}.$$

\begin{thm}\label{qg1}
Let $q > 1$, $\nu > -1$, $t=\frac{\nu + 2}{q}$, $\beta > \max (\frac{\nu}{q}, \frac{\nu + 2}{q} - 1)$ and
$f \in A^\infty_{\frac{\nu + 2}{q}}(\mathbb C_+)$. Then $l_1 = l_2$ where
$$l_1 = {\rm dist}_{A^\infty_{\frac{\nu + 2}{q}}}(f, A^q_\nu),$$

\begin{equation}\label{ggeps}
l_2 = \inf \left\{ \epsilon > 0 : \int_{\mathbb C_+} \left( \int_{V_{\epsilon, t}(f)}
\frac{ (\Im w)^{\beta -t} dm_2(w)}{|\overline w - z|^{\beta + 2}} \right)^q (\Im z)^\nu dm_2(z) < \infty \right\}.
\end{equation}
\end{thm}

\begin{thm}\label{ql1}
Let $0<q\leq 1$, $\nu > -1$, $t=\frac{\nu +2}{q}$, $\beta > \frac{\nu + 2}{q} - 2$ and
$f \in A^\infty_{\frac{\nu + 2}{q}}(\mathbb C_+)$. Then $t_1 = t_2$ where
$$t_1 = {\rm dist}_{A^\infty_{\frac{\nu + 2}{q}}}(f, A^q_\nu),$$

\begin{equation}\label{gleps}
t_2 = \inf \left\{ \epsilon > 0 : \int_{\mathbb C_+} \left( \int_{V_{\epsilon, t}(f)}
\frac{ (\Im w)^{\beta -t} dm_2(w)}{|\overline w - z|^{\beta + 2}} \right)^q (\Im z)^\nu dm_2(z) < \infty \right\}.
\end{equation}
\end{thm}

We present now the proofs of the above theorems, since they significantly overlap a unified presentation is possible. Assume $l_1 < l_2$ or $t_1 < t_2$. Then there are $\epsilon > \epsilon_1 > 0$ and $f_{\epsilon_1} \in A^q_\nu (\mathbb C_+)$ such that $\| f - f_{\epsilon_1} \|_{A^\infty_{\frac{\nu +2}{q}}} \leq \epsilon_1$ and
$$ \int_{\mathbb C_+} \left( \int_{V_{\epsilon, t}(f)}
\frac{(\Im w)^{\beta - t} dm_2(w)}{|\overline w - z|^{\beta + 2}} \right)^q (\Im z)^\nu dm_2(z) = +\infty.$$
Since $\| f - f_{\epsilon_1} \|_{A^\infty_{\frac{\nu +2}{q}}} \leq \epsilon_1$, from the definition of the
set $V_{\epsilon, t}(f)$ we conclude that
$$(\epsilon - \epsilon_1) \chi_{V_{\epsilon, t}(f)}(z) (\Im z)^{-\frac{\nu + 2}{q}} \leq |f_{\epsilon_1}(z)|$$
and therefore we have
\begin{eqnarray*}
+\infty & = & \int_{\mathbb C_+} \left( \int_{\mathbb C_+}
\frac{\chi_{V_{\epsilon, t}(f)}(w) (\Im w)^{\beta -t} dm_2(w)}{|\overline w - z|^{\beta + 2}} \right)^q
(\Im z)^\nu dm_2(z) \cr
& \leq & \int_{\mathbb C_+} \left( \int_{\mathbb C_+}
\frac{|f_{\epsilon_1}(w)| (\Im w)^\beta}{|\overline w - z|^{2+\beta}} dm_2(w)\right)^q (\Im z)^\nu dm_2(z) = M. \cr
\end{eqnarray*}
Note that this estimate is valid for $0<q<\infty$ and therefore works for both of the above theorems. In both cases we
are going to arrive at the contradiction by proving $M < +\infty$. Let us first consider the case $q>1$. Using
H\"older's inequality and Lemma \ref{lines} (with $\alpha = 0$) we obtain
\begin{eqnarray*}
I(z) & = & \left( \int_{\mathbb C_+} \frac{|f_{\epsilon_1}(w)|(\Im w)^\beta}{|\overline w - z|^{2+\beta}}
dm_2(w) \right)^q \cr
& \leq & \int_{\mathbb C_+} \frac{|f_{\epsilon_1}(w)|^q (\Im w)^{\beta q}}{|\overline w - z|^{\beta q - \epsilon q
+2}} dm_2(w) \cdot \left( \int_{\mathbb C_+} \frac{dm_2(w)}{|\overline w -z|^{\epsilon p + 2}} \right)^{q/p} \cr
& \leq & C \int_{\mathbb C_+} \frac{|f_{\epsilon_1}(w)|^q (\Im w)^{\beta q}}{|\overline w - z|^{\beta q - \epsilon q
+2}} dm_2(w) (\Im z)^{-\epsilon q}\cr
\end{eqnarray*}
for $w \in \mathbb C_+$, $\epsilon > 0$. Using this estimate, Fubini's theorem and Lemma \ref{lines} we obtain
\begin{eqnarray*}
M & \leq C & \int_{\mathbb C_+} |f_{\epsilon_1}(w)|^q (\Im w)^{\beta q} \left( \int_{\mathbb C_+}
\frac{(\Im z)^{\nu - \epsilon q}}{|\overline w - z|^{\beta q - \epsilon q + 2}} dm_2(z) \right) dm_2(w)\cr
& \leq & C \int_{\mathbb C_+} |f_{\epsilon_1}(w)|^q (\Im w)^\nu dm_2(w) < \infty \cr
\end{eqnarray*}
where $\epsilon > 0$ is small enough so that $\nu - \epsilon q > -1$.

Now let us turn to the case $q \leq 1$. We have, using Lemma \ref{lcub} and Lemma \ref{lines},
\begin{eqnarray*}
I(z) & = & \left( \int_{\mathbb C_+} \frac{|f_{\epsilon_1}(w)|(\Im w)^\beta}{|\overline w - z|^{2+\beta}}
dm_2(w) \right)^q \cr
& = & \left( \sum_k \int_{\Delta_k} \frac{|f_{\epsilon_1}(w)|(\Im w)^\beta}{|\overline w - z|^{2+\beta}}
dm_2(w) \right)^q \cr
& \leq & C \sum_{k=1}^\infty \max_{w \in \Delta_k} |f_{\epsilon_1}(w)|^q (m_2(\Delta_k))^q
\frac{(\Im w_k)^{\beta q}}{|\overline w - z|^{(2+\beta)q}} \cr
& \leq & C \int_{\mathbb C_+} \frac{|f_{\epsilon_1}(w)|^q (\Im w)^{\beta q + 2q - 2} dm_2(w)}{|\overline w - z|^{(2+
\beta)q}},\cr
\end{eqnarray*}
in the last estimate we used finite overlapping property of the family of enlarged cubes. Now we get $M < \infty$ as
in the case $q>1$, namely by applying Funini's theorem and integral estimate from Lemma \ref{lines}.

The reverse inequalities $l_1 \leq l_2$ and $t_1 \leq t_2$ can be proved simultaneously. We
fix $\epsilon > 0$ such that the integrals in (\ref{ggeps}), respectively (\ref{gleps}), are finite (if there
are no such $\epsilon > 0$, the inequality is trivial) and use integral representation from Lemma \ref{linre}. This integral representation is valid for $p = \infty$, $\alpha > 0$ and $\beta$ sufficiently large.
\begin{eqnarray*}
f(z) & = & \frac{\beta +1}{\pi}\left( \int_{\mathbb C_+ \setminus V_{\epsilon, t}(f)}
\frac{f(w) (\Im w)^\beta dm_2(w)}{(\overline w - z)^{2+\beta}} + \int_{V_{\epsilon, t}(f)}
\frac{f(w) (\Im w)^\beta dm_2(w)}{(\overline w - z)^{2+\beta}} \right) \cr
& = & \frac{\beta +1}{\pi}(f_1(z) + f_2(z)) \cr
\end{eqnarray*}
where $\beta$ is sufficiently large. The estimate for $f_1(z)$ is immediate, using Lemma \ref{lines}:
$$|f_1(z)| \leq \epsilon \int_{\mathbb C_+} \frac{(\Im w)^{\beta - t} dm_2(w)}{|\overline w - z|^{2+\beta}}
\leq C\epsilon (\Im z)^{-t},$$
and this means $f_1 \in A^\infty_{\frac{\nu + 2}{q}}$ with $\| f_1 \|_{A^\infty_t} \leq C \epsilon$. Next, by the
choice of $\epsilon > 0$ we have:
\begin{eqnarray*}
\| f_2 \|_{A^q_\nu} & \leq & C \| f \|_{A^\infty_t} \int_{\mathbb C_+} \left( \int_{V_{\epsilon, t}(f)}
\frac{(\Im w)^{\beta - t} dm_2(w)}{|\overline w - z|^{2+\beta}} \right)^q (\Im z)^\nu dm_2(z) \cr
& \leq &  C \| f \|_{A^\infty_t}.
\end{eqnarray*}
This completes the proofs of both theorems.

We end this section with an observation that the upper halfplane is the simplest tube domain. Analysis on tube domains
over symmetric cones is an active area of research (see \cite{BeBe}, \cite{BB} and \cite{BBG}) and availability of integral representations points to possibility of extending our results to this more general setting.

\section{Distance problems in $A^q_s(\Omega)$ spaces}

We now consider Bergman type spaces in $D \subset \mathbb C^n$, where $D$ is smoothly bounded relatively compact
strictly pseudoconvex domain, providing also sharp results in this case. Our proofs are heavily based on the estimates from \cite{Be}, where more general situation was considered.

Since $|f(z)|^p$ is subharmonic (even plurisubharmonic) for a holomorphic $f$, we have $A^p_s(D) \subset A^\infty_t(D)$ for $0<p<\infty$, $sp > n$ and $t = s$. Also, $A^p_s(D) \subset A^1_s(D)$ for $0 < p \leq 1$ and $A^p_s(D) \subset A^1_t(D)$ for $p > 1$ and $t$ sufficiently large. Therefore we have an integral representation
\begin{equation}\label{inre}
f(z) = \int_{D} f(\xi)K(z,\xi)\delta^t(\xi)dV(\xi),
\end{equation}
where $K(z,\xi)$ is a kernel of type $n+t+1$, that is a measurable function on $D \times D$ such that
$|K(z,\xi)| \leq C |\tilde \Phi(z,\xi)|^{-(n+1+t)}$, where $\tilde \Phi(z, \xi)$ is so called Henkin-Ramirez function
for $D$. From now on we work with a fixed Henkin-Ramirez function $\tilde \Phi$ and a fixed kerenel $K$ of type
$n+t+1$. We are going to use the following results from \cite{Be}.

\begin{lem}\label{be1}{\rm (\cite{Be}, Corollary 5.3.)}
If $r>0$, $0<p\leq 1$, $s>-1$, $p(s+n+1)>n$ and $f \in H(D)$, then we have
$$\left( \int_{D} |f(\xi)| |\tilde \Phi(z,\xi)|^r \delta^s(\xi) dV(\xi) \right)^p \leq C \int_{D}
|f(\xi)|^p |\tilde \Phi(z,\xi)|^{rp} \delta^{p(s+n+1)-(n+1)}(\xi) dV(\xi).$$
\end{lem}

\begin{lem}\label{be2}{\rm (\cite{Be}, Corollary 3.9.)}
Assume $T(z,\xi)$ is a kernel of type $\beta$, and $\sigma > 0$ satisfies
$\sigma + n -\beta < 0$. Then we have
$$\int_{D} T(z, \xi)\delta^{\sigma - 1}(z) dV(z) \leq C\delta^{\sigma + n - \beta}(\xi).$$
\end{lem}

A natural problem is to estimate ${\rm dist}_{A^\infty_s(D)}(f, A^q_s(D))$ where
$0<q<\infty$, $sq>n$ and $f \in A^\infty_s(D)$. We give sharp estimates below, treating cases $0<q\leq 1$ and
$q>1$ separately.

\begin{thm}
Let $0<q\leq 1$, $sq>n$, $f \in A^\infty_s(D)$ and $t > s$ is sufficiently large. Then $\omega_1 = \omega_2$
where
$$\omega_1 = {\rm dist}_{A^\infty_s(D)}(f, A^q_s(D)),$$

$$\omega_2 = \inf \left\{ \epsilon > 0 : \int_{D} \left( \int_{\Omega_{\epsilon, s}} |K(z,\xi)| \delta^{t-s}
(\xi) dV(\xi) \right)^q \delta^{sq-n-1}(z)dV(z) < \infty \right\},$$
where $K(z,\xi)$ is the above kerenel of type $n+t+1$ and
$$\Omega_{\epsilon, s} = \{ z \in D : |f(z)|\delta^s(z) \geq \epsilon \}.$$
\end{thm}

{\it Proof.} Let us prove that $\omega_1 \leq \omega_2$. We fix $\epsilon > 0$ such that the above integral
is finite and use (\ref{inre}):
$$f(z) = \int_{D \setminus \Omega_{\epsilon, s}} f(\xi)K(z,\xi) dV(\xi) + \int_{\Omega_{\epsilon, s}}
f(\xi)K(z,\xi) dV(\xi) = f_1(z) + f_2(z).$$
We estimate $f_1$:
\begin{eqnarray*}
|f_1(z)| & \leq  & C \epsilon \int_{D} |K(z,\xi)| \delta^{t-s} dV(\xi) \cr
& \leq & C \epsilon \int_{D} \frac{\delta^{t-s}(\xi) dV(\xi)}{|\tilde \Phi(z,\xi)|^{n+t+1}} \cr
& \leq & C \epsilon \delta^{-s}(z), \cr
\end{eqnarray*}
where the last estimate is contained in \cite{Be} (see p. 375). Next,
\begin{eqnarray*}
\| f_2 \|_{A^q_s}^q & = & \int_{D} |f_2(z)|^q \delta^{sq-n-1}(z)dV(z) \cr
& \leq & C \int_{D} \left(\int_{\Omega_{\epsilon, s}} |f(\xi)| K(z,\xi) \delta^t(\xi) dV(\xi) \right)^q
\delta^{sq-n-1}(z) dV(z)\cr
& \leq & C^\prime \| f \|^q_{A^\infty_s}.\cr
\end{eqnarray*}
Now we have
$${\rm dist}_{A^\infty_s(D)}(f, A^q_s(D)) \leq \| f - f_2 \|_{A^\infty_s(D)} =
\| f_1 \|_{A^\infty_s(D)} \leq C \epsilon.$$
Now assume that $\omega_1 < \omega_2$. Then there are  $\epsilon > \epsilon_1 > 0$ and $f_{\epsilon_1} \in A^q_s (D)$ such that $\| f - f_{\epsilon_1} \|_{A^\infty_s} \leq \epsilon_1$ and
$$I = \int_{D} \left( \int_{\Omega_{\epsilon, s}} |K(z,\xi)|\delta^{t-s}(\xi) dV(\xi) \right)^q
\delta^{sq-n-1}(z) dV(z) = \infty.$$
As in the case of the upper half-plane one uses $\| f - f_{\epsilon_1} \|_{A^\infty_s} \leq \epsilon_1$ to obtain
$$(\epsilon - \epsilon_1) \chi_{\Omega_{\epsilon, s}}(z) \delta^{-s}(z) \leq C|f_{\epsilon_1}(z)|.$$
Now the following chain of estimates leads to a contradiction:
\begin{eqnarray}\label{cont}
I & = & \int_{D} \left( \int_{D} \chi_{\Omega_{\epsilon, s}}(\xi) \delta^{t-s}(\xi)K(z,\xi) dV(\xi) \right)^q \delta^{sq-n-1}(z) dV(z) \cr
& \leq & C \int_{D} \left( \int_{D} |f_{\epsilon_1}(\xi)| \delta^t(\xi) K(z,\xi) dV(\xi) \right)^q
\delta^{sq-n-1}(z) dV(z) \cr
& \leq & C \int_{D} \left( \int_{D} |f_{\epsilon_1}(\xi)| \delta^t(\xi)
\frac{dV(\xi)}{|\tilde \Phi(z, \xi)|^{n+t+1}} \right)^q \delta^{sq-n-1}(z) dV(z)\cr
& \leq & C \int_{D} \int_{D} |f_{\epsilon_1}(\xi)|^q
\frac{\delta^{sq-n-1}(z) \delta^{q(t+n+1)-(n+1)}(\xi)}{|\tilde \Phi(z, \xi)|^{q(n+t+1)}} dV(z) dV(\xi) \cr
& \leq & C \int_{D} |f_{\epsilon_1}(\xi)|^q \delta^{sq-n-1}(\xi) dV(\xi) < \infty, \cr
\end{eqnarray}
where we used Lemma \ref{be1} and Lemma \ref{be2} with $\beta = q(n+t+1)$, $\sigma = sq-n$. $\Box$

Next theorem deals with the case $1 < q < \infty$.

\begin{thm}
Let $q>1$, $sq > n$, $t>s$, $t > \frac{s+n+1}{q}$ and $f \in A^\infty_s(D)$. Then $\omega_1 = \omega_2$ where
$$\omega_1 = {\rm dist}_{A^\infty_s(D)}(f, A^q_s(D)),$$

$$\omega_2 = \inf \left\{ \epsilon > 0 : \int_{D} \left( \int_{\Omega_{\epsilon, s}} |K(z,\xi)| \delta^{t-s}
(\xi) dV(\xi) \right)^q \delta^{sq-n-1}(z)dV(z) < \infty \right\}.$$
\end{thm}

{\it Proof.} An inspection of the proof of the previous theorem shows that it extends to this case also, provided
one can prove the estimate:
\begin{eqnarray*}
J & = & \int_{D} \left( \int_{D} |f_{\epsilon_1}(\xi)| \delta^t(\xi) K(z,\xi) dV(\xi) \right)^q
\delta^{sq-n-1}(z) dV(z) \cr
&\leq & C \int_{D} |f_{\epsilon_1}(\xi)|^q \delta^{sq-n-1}(\xi) dV(\xi) < \infty\cr
\end{eqnarray*}
where $q>1$. Using H\"older's inequality and Lemma \ref{be2}, with $\sigma = 1$ and
$\beta = n+1+p\epsilon$, we obtain
\begin{eqnarray*}
I(z) & = & \left( \int_{D} |f_{\epsilon_1}(\xi)| \delta^t(\xi) K(z,\xi) dV(\xi) \right)^q \cr
& \leq & \int_{D} \frac{|f_{\epsilon_1}(\xi)|^q \delta^{tq}(\xi) dV(\xi)}{|\tilde \Phi(z,\xi)|^{n+1+tq-
\epsilon q}} \cdot \left( \int_{D} \frac{dV(\xi)}{|\tilde \Phi(z, \xi)|^{n+1+p\epsilon}} \right)^{q/p}\cr
& \leq & C \int_{D} \frac{|f_{\epsilon_1}(\xi)|^q \delta^{tq}(\xi) dV(\xi)}{|\tilde \Phi(z,\xi)|^{n+1+tq-
\epsilon q}}\; \delta^{-q\epsilon}(z), \cr
\end{eqnarray*}
and this gives
\begin{eqnarray*}
J & \leq & C \int_{D} \int_{D} \frac{|f_{\epsilon_1}(\xi)|^q \delta^{tq}(\xi)\delta^{-q\epsilon +sq-n-1}(z)}{|\tilde \Phi(z,\xi)|^{n+1+tq-\epsilon q}} dV(z)\,dV(\xi) \cr
& \leq & C \int_{D} |f_{\epsilon_1}(\xi)|^q \delta^{sq-n-1}(\xi) dV(\xi) < \infty,\cr
\end{eqnarray*}
where we again used Lemma \ref{be2}, with $\beta = n+1+tq-\epsilon q$ and $\sigma = q(s-\epsilon)-n>0$. $\Box$

{\it Remark.} We note that most results of this paper and the previous one (\cite{SM}) on distances can be extended
to bounded symmetric domains $\Omega \subset \mathbb C^n$. Indeed, the methods of proofs are based on Bergman
representation formula, asymptotic properties of the Bergman kernel and Forelli-Rudin type estimates for integrals
$$\int_\Omega K^\alpha (w,w) |K(z,w)|^\beta dV(w), \;\; \alpha > 0, \; \beta > 0 \;\; z \in \Omega,$$
where $K(z,w)$ is a Bergman reproducing kernel for the weighted Bergman space $A^2_\alpha (\Omega)$. The relevant
estimates can be found in \cite{BeBe} and \cite{FK}.

\end{document}